\documentclass[12pt,a4paper]{article}
\usepackage[T1]{fontenc}
\usepackage{amsmath}
\usepackage{amssymb}
\usepackage{textcomp}
\usepackage[english]{babel}
\usepackage[all,dvips]{xy}
\usepackage[applemac]{inputenc}
\usepackage{xcolor}

\newtheorem{remk}{Remark}
 \newtheorem{hyp}{Assumption}

\newtheorem{thm}{Theorem}
\newtheorem{lem}[thm]{Lemma}
\newtheorem{prop}{Proposition}
\newtheorem{cor}[thm]{Corollary}

\newcommand\dem{\textbf{Proof: }}

\newcommand\rem{\textbf{Remark: }}
\newcommand\ex{\textbf{Example: }}

\title{Analytic linearization of a generalization of the semi-standard map: radius of convergence and Brjuno sum}
\author{C.Chavaudret, S.Marmi}
\date{}

\begin{document}

\maketitle

\textbf{Abstract:} One considers a system on $\mathbb{C}^2$ close to an invariant curve which can be viewed as a generalization of the semi-standard map to a trigonometric polynomial with many Fourier modes. The radius of convergence of an analytic linearization of the system around the invariant curve is bounded from below by $\exp(-\frac{2}{d}B(d\alpha)-C)$, where $C\geq 0$ does not depend on $\alpha$, $d\in \mathbb{N}^*$ and $\alpha$ is the frequency of the linear part. For a class of trigonometric polynomials, it is also bounded from above by a similar function. The error function is non decreasing with respect to the smallest coefficient of the trigonometric polynomial.

\section{Introduction}

Consider the following discrete dynamical system:

\begin{equation}\label{systemediscret}
\left\{ \begin{array}{l}
x_{n+1} = x_n+y_n+ A(x_n)\\
y_{n+1} = y_n+A(x_n)\\
\end{array}\right.
\end{equation}

\noindent where $A$ is a trigonometric polynomial with positive Fourier modes: $A(x) = \sum_{K=1}^N a_K e^{iKx}$, $N\in\mathbb{N}$, $a_K\in \mathbb{C}$. 
This can be viewed as a generalization of the semi-standard map, which is given by $A(x)=e^{ix}$. 
The semi-standard map was introduced as a model which is similar to, but easier to study than, the standard map given by $A(x)=\sin x$, where the coexistence of both positive and negative Fourier modes makes the linearization problem more difficult to solve. 


With the change of variables $z_n=e^{ix_n}, \lambda_n=e^{iy_n}$, system \eqref{systemediscret} is conjugated to the system

\begin{equation}\label{systemeinitial}
\left\{ \begin{array}{l}
z_{n+1} = \lambda_n z_n \prod_{K=1}^N e^{i a_K z_n^K}\\
\lambda_{n+1} = \lambda_n \prod_{K=1}^N e^{i a_K z_n^K}
\end{array}\right.
\end{equation}

\noindent for which the set $\{0\}\times S^1$ is invariant.
Let $F(\lambda,z)=(\lambda z \prod_{K=1}^N e^{ia_Kz^K},\lambda  \prod_{K=1}^N e^{ia_Kz^K})$ so the system \eqref{systemeinitial} can be written $(z_{n+1},\lambda_{n+1})=
F(z_n,\lambda_n)$. One looks for an analytic linearization of the system, that is to say, an analytic map $H(z,\lambda)= (h(z,\lambda),h_2(z,\lambda))$ such that $F\circ H = H\circ R$ where $R(z,\lambda)= (\lambda z,\lambda)$. 
If there is such a linearization and if it is analytic, then the invariant curves of the rotation $R$ are smoothly preserved. 
Then one has a family of invariant closed curves in $\mathbb{C}^2$, corresponding to the numbers $\lambda$ with modulus 1 and $z$ in a neighbourhood of $0$ where the linearization is analytic. Now if $\lambda$ is rational, 
there is no way of finding a linearization which is analytic in $z$. 
Thus we will have to assign a value to the parameter $\lambda$, with $|\lambda|=1$ and $arg(\lambda)\in\mathbb{R}\setminus\mathbb{Q}$, to construct a linearization which is close to the identity and analytic; its radius of convergence will depend on the arithmetical properties of $\lambda$. 

\bigskip
Davie \cite{D94}Êand Marmi \cite{M90} proved that concerning the semi-standard map, the radius of convergence $\rho(\alpha)$ of the linearization is bounded as follows:

$$\exp(-2B(\alpha) -C)\leq \rho(\alpha)\leq \exp(-2B(\alpha)+C')$$

\noindent where $C>0,C'>0
$ do not depend on the complex argument $\alpha$ of $\lambda$ and where $B(\alpha)$ is the Brjuno sum of $\alpha$. In particular, if $B(\alpha)$ diverges, then there is no analytic linearization around $0$. This can be reformulated as stating that the error function $\alpha\mapsto 2B(\alpha)+\ln \rho(\alpha)$ is bounded.

\bigskip
\noindent
After the numerical evidence in \cite{MS92}, a similar result about the standard map, in the perturbative case, was proved in \cite{BG01} and \cite{BG02}. 
However, concerning the semi-standard map and in the present paper,
the strong assumption of Fourier modes being only positive makes it possible to remove the perturbative assumption.

\bigskip
\noindent
The Brjuno sum was first introduced in \cite{Br72} to give a sufficient condition to the convergence of the linearization for analytic vector fields around a fixed point. 

 \bigskip
\noindent
 Yoccoz proved in \cite{Yo88} that the Brjuno condition (i.e the convergence of the Brjuno function, or equivalently of the Brjuno sum) is necessary and sufficient to the analytic linearization of the quadratic polynomial and of germs of diffeomorphisms of $(\mathbb{C},0)$. This resulted in the study of the error function $\Phi+\ln r$, where $\Phi$ is the Brjuno function and $r$ is the radius of convergence of the linearization for the quadratic polynomial. It was conjectured in \cite{MMY97} that this function is $1/2$-Hlder and Buff and Chritat showed first that it is bounded in \cite{BC04}, then that it is continuous in \cite{BC06}. Cheraghi-Chritat then proved that a restriction of this error function is $1/2$-Hlder.  
 
 \bigskip
\noindent
Brjuno's lower bound on the convergence radius for linearization of analytic vector fields was improved in \cite{GM10}. In \cite{S00}, Stolovitch replaced the arithmetical condition on the spectrum of the linear part by a condition of algebraic nature.
The Brjuno sum was also proved to play a role in other analytic linearization problems, as for instance linearization of vector fields around an invariant torus (see \cite{Au13} and \cite{C16}), or reducibility of quasiperiodic cocycles (\cite{ChM12}).

\bigskip
\noindent The optimality of the Brjuno condition was also studied in other linearization problems. Carletti-Marmi proved in \cite{CaM00} that the Brjuno condition is also necessary to linearize analytically a germ of diffeomorphism around a fixed point, and generalize it to Gevrey classes.  However the continuity of the analogue of the error function for linearization problems more general than the quadratic polynomial remains open up to now.


\bigskip
\noindent
Our main result is stated in the following two theorems:

\begin{thm}
Let $\rho$ be the radius of convergence of the linearization of the system \eqref{systemeinitial}. Let $d$ be the greatest common divisor of the indices of the Fourier modes of $A$. 
There exists $C\geq 0$, which is a non decreasing function of the greatest coefficient of $A$
and which does not depend on $\alpha$,
such that

\begin{equation}\label{minoration}
\rho\geq \exp(-
\frac{2}{d}B(d\alpha)-C)
\end{equation}


\end{thm}

\begin{thm}\label{thm2}
Assume that the coefficients of $A$ satisfy the following assumption: there exists $\theta\in\mathbb{R}$ such that for all $k=1,\dots, N$ with $a_k\neq 0$, the complex argument of the number $a_k$ is $k\theta+\frac{\pi}{2}$. 

\noindent
Let $\rho$ be the radius of convergence of the linearization of the system \eqref{systemeinitial}. Let $d$ be the greatest common divisor of the indices of the Fourier modes of $A$. There exists $C'\geq 0$, which is a non decreasing function of the smallest coefficient of $A$ and does not depend on $\alpha$
such that

\begin{equation}\label{majoration}
\rho\leq  \exp(- \frac{2}{d}B(d\alpha) +C')
\end{equation}

\end{thm}

\noindent
In particular, if $B(d\alpha)$ diverges, then there is no analytic linearization.

\bigskip
\noindent \rem Let $\kappa_0$ be the smallest integer such that $a_{\kappa_0}\neq 0$. 
The assumption in Theorem \ref{thm2} says that the argument of every non zero $a_k$ is the following function of $a_{\kappa_0}$:

$$Arg(a_k)=\frac{k}{\kappa_0} Arg(a_{\kappa_0}) +\frac{\pi}{2}(1-\frac{k}{\kappa_0})$$

\noindent This assumption is in the spirit of Cremer's couterexample of non-linearizable germs, except that the arguments are defined from the beginning instead of recursively.


\bigskip
\noindent
The main result is obtained by a direct analysis of the coefficients of the formal linearization (which always exists, and is unique if one requires it to be formally close to the identity); as in \cite{M90} and \cite{D94},
it appears that there is a link between those coefficients and the Brjuno sum of $d\alpha$. 
To get an upper bound of the radius of convergence, that is to say, a lower bound on the coefficients of the linearization, one uses a strong assumption on the complex arguments of the 
coefficients of $A$, in order to be able to bound the sum from below by just one of its terms. 
However the lower bound on the radius of convergence does not use this assumption.

\section{Notations}
Let $x\in\mathbb{R}$, one denotes by $||x||_\mathbb{Z}$ the distance between $x$ and the closest integer: $||x||_\mathbb{Z}=\min_{p\in\mathbb{Z}}|x-p|$.

\bigskip
\noindent
Considering the trigonometric polynomial $A(x)$ defined at the beginning, denote by $\kappa_0 <\kappa_1<\dots <\kappa_N$ all indices of Fourier modes of $A$: thus for all $0\leq i\leq N$, $a_{\kappa_i}\neq 0$, and if $ \forall i=0,\dots,N$, $K\neq \kappa_i$, then $a_K=0$. Also denote by $\kappa_{min}$ (resp. $\kappa_{max}$) the number $\kappa_i$ minimizing (resp. maximizing) $\{|a_{\kappa_i}|,i=0,\dots, N\}$.

\bigskip
\noindent
Denote by $\mathcal{M}\subset \mathbb{N}$ the additive semi-group generated by $\{\kappa_0,\dots, \kappa_N\}$:

\begin{equation}\label{def-semigp}
\mathcal{M}=\{ p_0\kappa_0+\dots+p_N\kappa_N>0, p_0\in\mathbb{N},\dots, p_N\in\mathbb{N}\}
\end{equation}

\section{Analysis of the linearization}

With a reasoning similar to the one in \cite{BMS00} (which is reproduced and adapted to the present model in the appendix), 
one proves that if $H$ linearizes the system, then, up to a multiplicative constant which will not change the radius of convergence, the first component of the linearization is a function
$h(z,\lambda)=ize^{\Phi_\lambda(z)}$ satisfying (mod $2i\pi$):

\begin{equation}\label{equation-Phi}\sum_{K=1}^N ia_K(iz)^K e^{K\Phi_\lambda(z)} = \Phi_\lambda (\lambda^{-1} z) +\Phi_\lambda(\lambda z) -2\Phi_\lambda(z)
\end{equation}

\noindent
Let $\lambda\in\mathbb{C}$ with modulus 1 and $\alpha=\frac{Arg(\lambda)}{2\pi}$.
Denoting
$\Phi_\lambda(z) = \sum_{l\geq 1} \phi_l z^l$ and

\begin{equation}\label{diviseurs}d_{l,\lambda}=\lambda^{\frac{l}{2} } -\lambda^{\frac{-l}{2} } =2i\sin(\pi l\alpha)
\end{equation}

\noindent one gets the following:

\begin{equation}
\sum_{K\in\{\kappa_0,\dots,\kappa_N\}} i^{K+1} a_K z^K\left( 1+ \sum_{p\geq 1} \frac{1}{p!} (\sum_{j\geq 1} \phi_j z^j )^p\right)^K
=\sum_{l\geq 1} \phi_l z^l d_{l,\lambda}^2
\end{equation}

\noindent
In what follows, $\lambda$ being fixed with modulus 1 and argument $\alpha$ such that $\frac{\alpha}{2\pi}\in \mathbb{R}\setminus\mathbb{Q}$, we will denote for all $l\in\mathbb{N}\setminus\{0\}$,

$$D_l:=-4\sin^2(\pi l\alpha)<0$$

\noindent Then $|D_l|^{-1} \geq \frac{1}{4}$ for all $l\geq 1$.
Let $l\geq 1$, then $\phi_l=0$ if $l<\kappa_0$ and if $l\geq \kappa_0$,

\begin{equation}\label{rec}
\phi_l = \frac{1}{D_l}\sum_{\substack{K\in\\
\{\kappa_0,\dots,\\
\kappa_N\},\\
K\leq l}} i^{K+1} a_K[\delta_{K,l}+\sum_{\substack{m=1,\dots, K\\
m\leq l-K}} C^m_K \sum_{p_1=1}^{l-K} \dots \sum_{p_m= 1}^{l-K}\frac{1}{p_1!\dots p_m!}\sum_{
\substack{j_1^1,\dots , j^1_{p_1},\\
\vdots\\
j^m_1,\dots, j^m_{p_m},\\
j_1^1+\dots + j^1_{p_1} +\dots\\
 + j^m_1 +\dots 
+ j^m_{p_m} = l-K}} \phi_{j^1_1}\dots \phi_{j^1_{p_1}}\dots \phi_{j^m_1}\dots \phi_{j^m_{p_m}}
]
\end{equation} 

\bigskip
\noindent Another recurrence relation is required:
let $\Psi_\lambda(z)=\sum_{K\in \{\kappa_0,\dots,\kappa_N\}}\Psi_K(z), \ \Psi_K(z)=i^{K+1}a_Kz^Ke^{K\Phi_\lambda(z)}$ and denote by $\sum_{k\geq 1} \psi_kz^k$ the Taylor expansion of $\Psi_\lambda$ and by $\sum_l \psi_{K,l}z^l$ the Taylor expansion of $\Psi_K$. Then 

\begin{equation}\label{equation-PsiK}\Psi_\lambda(z) = \Phi_\lambda (\lambda^{-1} z) +\Phi_\lambda(\lambda z) -2\Phi_\lambda(z)
\end{equation}

\noindent which implies that for all $l\geq \kappa_0$,

\begin{equation}\label{recPhi}
\psi_l=D_l\phi_l
\end{equation}

\noindent Moreover, derivating $\Psi_K$,
one sees that $\psi_{K,l}=0$ for all $l<K$, $\psi_{K,K}=i^{K+1}a_K$, and
for all $n\geq K+1$,

\begin{equation}\label{recPsi}(n-K)\psi_{K,n} = K\sum_{k=1}^{n-1} k\phi_k \psi_{K,n-k}
\end{equation}

%

\bigskip
\rem One sees that $\phi_l\neq 0 \Rightarrow l\in \mathcal{M}$. This can be shown by recurrence: for $l=\kappa_0$, the property holds. Assume it holds for all $\kappa_0\leq l'\leq l-1$, for a fixed $l>\kappa_0$. Assume that $\phi_l\neq 0$. If $l\in\{\kappa_0,\dots,\kappa_N\}$, then $l\in\mathcal{M}$; otherwise, there exists $K\in\{\kappa_0,\dots,\kappa_N\}$, $m\leq l-K$, $p_1\leq l-K, \dots, p_m\leq l-K$ and non vanishing $\phi_{j_1^1},\dots, \phi_{j_{p_m}^m}$ such that $j_1^1+\dots + j_{p_m}^m = l-K$. By recurrence assumption, $j_1^1,\dots, j_{p_m}^m\in\mathcal{M}$ therefore $l\in\mathcal{M}$. 

\bigskip
\noindent A similar fact holds for every $\Psi_K$: if $\psi_{K,l}\neq 0, l\geq K,$ then $l\in\mathcal{M}$. Indeed 
$\psi_{K,l}=0$ if $l< K$ 
so the recurrence property holds for $l\leq K$. Assume the property holds up to a fixed $l\geq K$. If $\psi_{K,l+1}\neq 0$ then by \eqref{recPsi} there exist $k\in\{1,\dots,n-1\}$ such that $\phi_k\neq 0$ and $\psi_{K,l+1-k}\neq 0$ so $k\in\mathcal{M}$ and $ l+1-k\in\mathcal{M}$, therefore $l+1\in\mathcal{M}$.

\bigskip
\ex If $A$ is a monomial of order $K$, then $\Phi_\lambda(z)$ only has coefficients indexed by multiples of $K$. Then there exists a function $\Xi_\lambda:\mathbb{C}\rightarrow \mathbb{C}$ such that for all $z$, $\Phi_\lambda(z)=\Xi_\lambda(z^K)$.


\bigskip
\noindent The following lemma will be used to bound the coefficients of $\Phi$ from below by one of the terms of the sum determining them.

\begin{lem}\label{argument}
If there exists $\theta\in\mathbb{R}$ such that for all $K\in\{\kappa_0,\dots,\kappa_N\}$, the complex number $a_K $ has argument $K\theta+\frac{\pi}{2}$, 
then for all $l\in \mathcal{M}$, 
$\phi_l$ has argument $l(\theta+\frac{\pi}{2})$ and for all $K\in\{\kappa_0,\dots,\kappa_N\}$, if $\psi_{K,l}\neq 0$ then $\psi_{K,l}$ has argument $l(\theta+\frac{\pi}{2})+\pi$.
\end{lem}

\dem First one proves the part of the statement concerning $\phi_l$:
\begin{enumerate}
\item \label{1} If $l=\kappa_0$,
then \eqref{rec} implies that $\phi_l=\frac{1}{D_{\kappa_0}} i^{\kappa_0+1}a_{\kappa_0}$ which is, by assumption, the product of a number of modulus
$\kappa_0 (\theta+\frac{\pi}{2})$ with $\frac{-1}{D_{\kappa_0}}$, the latter being real and positive. 
\item \label{2} 
Let $l\geq \kappa_0+1$. Assume that for all $l'\in\mathcal{M}$ such that $l'\leq l-1$, $\phi_{l'}$ has argument $l'(\theta+\frac{\pi}{2})$. If $l$ differs from all $\kappa_i$, then by \eqref{rec}, $\phi_l$ is the sum of terms which are the product of a number $\frac{i^{K+1}a_K}{D_l}$, which by assumption has argument $K(\theta+\frac{\pi}{2})$, with a number which by recurrence assumption has 
argument $(l-K)(\theta+\frac{\pi}{2})$. 
If $l$ is one of the $\kappa_i$, then one has to add the term $\frac{1}{D_l}i^{l+1}a_l$, which has argument $l(\theta+\frac{\pi}{2})$.
\end{enumerate}

As for the statement concerning $\psi_{K,l}$, since $\psi_{K,K}=i^{K+1}a_K$ which has argument $K(\theta+\frac{\pi}{2})+\pi$, the property holds for $l\leq K$. Assume this property holds up to a fixed $l\geq K$. By equation \eqref{recPsi}, if $\psi_{K,l+1}\neq 0$ then $\psi_{K,l+1}$ is a sum of terms with argument $k(\theta+\frac{\pi}{2}) + (l+1-k)(\theta+\frac{\pi}{2})+\pi$ for $k\in\{1,\dots, n-1\}$, therefore it has argument $(l+1)(\theta+\frac{\pi}{2})+\pi$.
$\Box$

\bigskip
\rem 
For instance, $A(x)$ satisfies the assumption of Lemma \ref{argument} if it only has coefficients $a_K\in i\mathbb{R}^+$ (then $\theta= 0$), or such that $i^{k+1}a_K\in\mathbb{R}^-$ (then $\theta=-\frac{\pi}{2}$), or if $A$ is a monomial (without restriction on $\theta$).

\subsection{The semigroup $\mathcal{M}$}

We shall also need the following lemma on the set $\mathcal{M}$:

\begin{lem}\label{semigroupe}Let $d$ be the greatest common divisor of $\kappa_0,\dots,\kappa_n$. There exists $N_\mathcal{M}$ such that for all integer $m\geq N_\mathcal{M}$, if $m$ is a multiple of $d$, then $m\in\mathcal{M}$. If $\kappa_0=1$ then $N_\mathcal{M}=1$.
\end{lem}

\dem Let $I_0=[\kappa_0,\kappa_0+\dots+\kappa_N[\cap \mathbb{N}$ and for $p\geq 1$, let

$$I_p=[p(\kappa_0+\dots+\kappa_N),(p+1)(\kappa_0+\dots+\kappa_N)[\cap\mathbb{N}$$

\noindent Then
$(I_p)_{p\geq 0}$ is a partition of $[\kappa_0,+\infty[\cap \mathbb{N}$. Moreover, if $jd\in I_p\cap\mathcal{M}$ for some $p\geq 0$, then $jd+\kappa_0+\dots +\kappa_N$, which is also a multiple of $d$, belongs to $I_{p+1}\cap \mathcal{M}$; therefore, in order to prove this lemma, it is sufficient to prove that for all $jd\in I_0\setminus\mathcal{M}$, there exists $p\geq 1$ such that $jd+p(\kappa_0+\dots +\kappa_N)\in\mathcal{M}$. After a finite number of steps, one obtains an integer $P$ such that all multiples of $d$ belonging to $I_P$ are also elements of $\mathcal{M}$. Translating by $\kappa_0+\dots + \kappa_N$, one will deduce that all multiples of $d$ greater than $p(\kappa_0+\dots +\kappa_N)$ are also elements of $\mathcal{M}$.

\bigskip
Now, let $jd\in I_0\setminus\mathcal{M}$ if it exists (otherwise the proof is finished). By Bezout's theorem, there are relative integers $b_0,\dots , b_N$, at least one of which is positive, such that $b_0 \kappa_0+\dots + b_N\kappa_N = jd$. Sorting the coefficients $b_i$ by sign, one infers that $jd-\sum_{i/ \ b_i\leq 0} b_i\kappa_i$ is a linear combination of the $\kappa_i$ with non negative coefficients, at least one of which is non zero, therefore $jd-\sum_{i/ \ b_i\leq 0}b_i\kappa_i\in\mathcal{M}$. Therefore $jd +\max(-b_0,\dots , -b_N)(\kappa_0+\dots+\kappa_N)\in\mathcal{M}$. $\Box$

\bigskip
\noindent Lemma \ref{semigroupe} has the following corollary:

\begin{cor}\label{corsemigp}
There exists $N_\mathcal{M}$ such that for all $a,b\in\mathcal{M}$, if $a-b\geq N_\mathcal{M}$, then $a-b\in\mathcal{M}$.
\end{cor}

\dem Let $d$ be the greatest common divisor of $\kappa_0,\dots,\kappa_N$. Since every element of $\mathcal{M}$ is a multiple of $d$, if $a,b\in\mathcal{M}$, then $a-b$ is a multiple of $d$. If moreover $a-b\geq N_\mathcal{M}$, where $N_\mathcal{M}$ was defined in Lemma \ref{semigroupe}, then $a-b\in\mathcal{M}$. $\Box$

\section{The Brjuno sum and the Brjuno function}

Let $d$ be the greatest common divisor of $\kappa_0,\dots,\kappa_N$.
Let us consider the continued fraction expansion of $d\alpha$.

\bigskip
\noindent
\textbf{Notations:} 
Let $(q_k)$ be the sequence of the denominators of the approximants of $d\alpha$.
Recall the well-known recurrence relation: for all $j\geq 0$,

\begin{equation}\label{suite-denom}
q_{j+2}=a_{j+2} q_{j+1} + q_j
\end{equation}

\noindent where $(a_j)$ is the sequence of integers given by the continued fraction expansion.


%
%
%

%
%

\bigskip
\noindent
The following lemmas are given in order to relate the Brjuno sum with the small divisors of our linearization problem.

\begin{lem}\label{grand-div}
For all $k\geq 1$, there is

$$\frac{1}{2q_{k+1}}\leq |D_{dq_k}|^{\frac{1}{2}} \leq \frac{3}{q_{k+1}}$$

\end{lem}

\dem For all $l\in \mathbb{Z}$,

$$|D_l| = 4|\sin(\pi l\alpha)|^2 = 4 [(||l\alpha||_\mathbb{Z} + R(|| l\alpha||_\mathbb{Z}))]^2$$

\noindent (where $R$ is the remainder in the Taylor-Lagrange formula)
whence

\begin{equation}\label{encadre-Dl}
|| l\alpha || \leq |D_l|^{\frac{1}{2}}\leq 3 || l\alpha ||
\end{equation}

\noindent
Now for all $k\geq 1$,

$$ ||q_kd\alpha||_\mathbb{Z}=\min_{p\in\mathbb{Z}}|q_kd\alpha-p|=q_k\min_{p\in\mathbb{Z}}|d\alpha-\frac{p}{q_k}|=q_k|d\alpha-\frac{p_k}{q_k}|$$

\noindent
Since

\begin{equation}\label{approx-rationnels} \frac{1}{2q_kq_{k+1}}\leq |d\alpha - \frac{p_k}{q_k}|\leq \frac{1}{q_kq_{k+1}}
\end{equation}

\noindent (see for instance \cite{MMY97} remark 1.6), there is

$$\frac{1}{2q_{k+1}}\leq |D_{dq_k}|^{\frac{1}{2}} \leq \frac{3}{q_{k+1}}.\ \Box$$

\rem This cannot be extended to generalized continued fractions since the inequality \eqref{approx-rationnels} does not hold anymore for generalized continued fractions.

\bigskip
\noindent The following lemmas come from \cite{D94}. They concern the structure of the set of small divisors and will be used in the lower bound on the radius of convergence.
We recall them here and apply the second to our setting.

\begin{lem}\label{Davie1}(\cite{D94}, lemma 2.2) 
Let $k\in \mathbb{N}$ and $n\in\mathbb{N},n\geq 1$, such that $d$ is a divisor of $n$. If $|D_n|^{\frac{1}{2}} <\frac{1}{q_k}$, then $\frac{n}{d}\geq q_k$ and either $q_k$ divides $\frac{n}{d}$ or $\frac{n}{d}\geq \frac{q_{k+1}}{4}$.
\end{lem}

\begin{lem}\label{Davie2}(\cite{D94}, lemma 2.3)
For all $k\geq 0, n\geq 1$, let $A_k(n)=\{dq_k\leq  j\leq n,\ d|j,\ \frac{1}{6q_{k+1}}\leq |D_{j}|^{\frac{1}{2}}<\frac{1}{6q_k}\}$. Let $E=\max(dq_k,\frac{dq_{k+1}}{4})$. Then 
there is a function $g_k:\mathbb{N}\rightarrow \mathbb{R}^+$ such that:

\begin{itemize}
\item $g_k(n)\leq (1+\frac{2dq_k}{E})\frac{n}{dq_k}$;
\item for all $n_1,n_2\in \mathbb{N}$, $g_k(n_1)+g_k(n_2)\leq g_k(n_1+n_2)$;
\item if $n\in A_k(n)$, then $g_k(n)\geq g_k(n-1)+1$.
\end{itemize}
\end{lem}

\rem This application of Davie's lemma is possible because the set $A_k$ satisfies: if $j_1=dj'_1<j_2=dj'_2\in A_k$, then either $dq_k$ divides $j_2-j_1$ or $j_2-j_1\geq \frac{dq_{k+1}}{4}$. Indeed, letting $p_1$ (resp. $p_2$) be the integer closest to $j'_1d\alpha$ (resp. $j'_2d\alpha$),

$$||(j_2-j_1)\alpha|| \leq |p_2-p_1 -(j'_2-j'_1)d\alpha|\leq ||j'_2d\alpha || +||j'_1d\alpha||\leq |D_{j_2}|^{\frac{1}{2}}+|D_{j_1}|^{\frac{1}{2}}< \frac{1}{3q_k}$$

\noindent (the last inequality comes from $j_1,j_2\in A_k$).
Therefore, by \eqref{encadre-Dl}, $|D_{j_2-j_1}| ^{\frac{1}{2}}<\frac{1}{q_k}$. Lemma \ref{Davie1} then implies that $q_k$ divides $j'_2-j'_1$ or $j'_2-j'_1\geq \frac{q_{k+1}}{4}$.

\begin{lem}\label{somme-convergente}
There exists $C_0,C_1>0$ such that

$$\sum_{l\geq 0}\frac{1}{q_{{l}}} \leq C_0$$

\noindent and

$$\sum_{l\geq 0}\frac{1}{q_{{l}}} \ln q_{{l}}\leq C_1
$$
\end{lem}

\dem Let us give a bound on $\sum_{l\geq 0}\frac{1}{q_{{l}}} \ln q_{{l}}$. Because of the recurrence relation \eqref{suite-denom}, the sequence $(q_{l})$ increases at least as fast as a Fibonacci sequence the first two terms of which are in $\{\kappa_0,\dots,\kappa_N\}$: denoting by $(f_k)$ the Fibonacci sequence with $f_0=f_1=1$, one recursively proves that

\begin{equation}\label{min-fibonacci}q_{k}\geq f_k
\end{equation}

\noindent 
Indeed, $q_{0}\geq 1$
and $q_{1}\geq 1$.
Assume that $q_{{k-1}} \geq 
f_{k-1} $ and $q_{{k}} \geq 
f_{k} $, then

$$q_{{k+1}} = a_{{k+1}}q_{{k+1}-1}+q_{{k+1}-2}\geq q_{k}+q_{{k-1}}\geq 
(f_k+f_{k-1})=
 f_{k+1}.$$

\noindent
Therefore, since the function $t\mapsto \frac{\ln t}{t}$ decreases on $[e,+\infty[$,

$$\sum_{l\geq 0}\frac{1}{q_{{l}}} \ln q_{{l}}\leq \ln 2+\sum_{l\geq 0} \frac{\ln f_l}{ f_l}\leq C_1
$$

\noindent where $C_1$ is a numerical constant.
From the inequality \eqref{min-fibonacci}, one also infers that

$$\sum_{l\geq 0}\frac{1}{q_{{l}}} \leq 
\sum_{l\geq 0} \frac{1}{f_k}\leq C_0
$$

\noindent where $C_0$ is a numerical constant. $\Box$

\subsection{Subsequence of fast increasing denominators}\label{subsequence}

Given the number $N_\mathcal{M}$ which was defined in Lemma \ref{semigroupe},
let $(n_k)$ the subsequence containing all indices such that

$$\left\{\begin{array}{l}
q_{n_0}\geq \max(N_\mathcal{M}+2,1+\kappa_{max}(2+\frac{1}{d}))>q_{n_0-1},\\
q_{n_{k}+1}\geq q_{n_k}^2+\zeta(\kappa_{max},d)q_{n_k} + \eta (\kappa_{max},d)\end{array}\right.$$

\noindent where $\zeta,\eta$ are given by

$$\zeta(\kappa_{max},d)=\frac{\kappa_{max}}{d}+3\kappa_{max}+2, \ \ \eta(\kappa_{max},d)=(\frac{\kappa_{max}}{d}+2\kappa_{max})(\kappa_{max}+1)$$

\begin{lem}\label{suite-extraite}
For all $k\geq 0$, it holds that
$q_{n_k}\geq \max(N_\mathcal{M}+2,\kappa_{max})^{2^k}$.

\end{lem}

\dem By definition, $q_{n_0}\geq \max(N_\mathcal{M}+2,\kappa_{max})$.

\noindent
Assume that $q_{n_k}\geq \max(N_\mathcal{M}+2,\kappa_{max})^{2^k}$, then $q_{n_{k+1}}\geq q_{n_k+1}\geq q_{n_k}^2\geq \max(N_\mathcal{M}+2,\kappa_{max})^{2^{k+1}} $. $\Box$

\begin{cor}\label{somme-suite-extraite}
For all $k\geq 0$, $\sum_{k\geq 0}\frac{1}{q_{n_k}}\leq \frac{1}{\max(N_\mathcal{M}+2,\kappa_{max})-1}$.
\end{cor}

\dem Indeed

\begin{equation}\begin{split}\sum_{k\geq 0}\frac{1}{q_{n_k}}&\leq \sum_{k\geq 0} \frac{1}{\max(N_\mathcal{M}+2,\kappa_{max})^{2^k}}\leq \sum_{k\geq 1} \frac{1}{\max(N_\mathcal{M}+2,\kappa_{max})^k}\\
&\leq 
 \frac{1}{\max(N_\mathcal{M}+2,\kappa_{max})-1}.\ \Box
\end{split}\end{equation}

\begin{lem}\label{somme-extraite}
There exists $C_2>0$ such that

$$B(d\alpha) -C_2\leq \sum_{l\geq 0} \frac{|\ln |D_{q_{n_l}}|\ |}{2q_{n_l}} \leq  B(d\alpha)+C_2' $$

\end{lem}

\dem Lemma \ref{grand-div} implies that

$$\sum_{l\geq 0}\frac{\ln q_{n_l+1}-\ln 3}{q_{n_l}} \leq \sum_{l\geq 0} \frac{|\ln |D_{q_{n_l}}|\ |}{2q_{n_l}} \leq \sum_{l\geq 0}\frac{\ln q_{n_l+1}}{q_{n_l}}+\sum_{l\geq 0} \frac{\ln 2}{q_{n_l}}
\leq B(d\alpha)+C_0\ln 2$$

\noindent \noindent (where we have also used Lemma \ref{somme-convergente}).
Now

$$\sum_{l\geq 0}\frac{\ln q_{n_l+1}}{q_{n_l}}=B(d\alpha) - \sum_{l\geq 0, q_l ^2+\zeta q_l +\eta> q_{l+1}}\frac{\ln q_{l+1}}{q_{l}}
$$

\noindent (where $\zeta, \eta$ were defined at the beginning of the section). Thus

\begin{equation}\begin{split}\sum_{l\geq 0}\frac{\ln q_{n_l+1}}{q_{n_l}}&\geq B(d\alpha) - \sum_{l\geq 0}\frac{\ln (q_l^2+\zeta q_l+\eta)}{q_{l}}\geq B(d\alpha)-\sum_{l\geq 0} 
\frac{\ln q_l + \ln (q_l+\zeta+\eta)}{q_l}\\
&\geq B(d\alpha) - 2\sum_{l\geq 0}\frac{\ln q_l}{q_l}-\sum_{l\geq 0} \frac{\ln (\zeta +\eta+1)}{q_l}
\end{split}\end{equation}

\noindent
By Lemma \ref{somme-convergente}, 

$$\sum_{l\geq 0}\frac{\ln (\zeta +\eta+1)}{q_{{l}}} \leq  C_0\ln (\zeta +\eta+1)
$$

\noindent
Therefore one can define $C_2 =  2C_1 +C_0(\ln 3+\ln (\zeta +\eta+1))$ and $C_2'=C_0\ln 2
$. $\Box$

\section{Recursively defined lower bound}

In this section we introduce a function ${F}$ which will be used in giving a lower bound on the coefficients of the linearization.
More precisely, we shall prove that

$$\limsup_{k\rightarrow +\infty} \frac{F(q_{n_k})}{q_{n_k}}-C\leq \limsup_{k\rightarrow +\infty} \frac{\ln |\phi_{q_{n_k}}|}{q_{n_k}}
 $$
 
 \noindent where $C$ 
 is a constant not depending on $\alpha$.

\noindent
Then we shall bound 
$\limsup_{k\rightarrow +\infty} \left(\frac{F\left(q_{n_k}\right)}{q_{n_k}}\right)$ 
by means of the Brjuno sum
in the Lemma
\ref{fonction-minorante-brjuno} 
below.



\bigskip
\noindent The function $F$ is recursively defined on the set $\{q_{n_k},k\geq 0\}$ as follows:

\begin{equation}\label{fonctionF}\begin{split}&F(q_{n_0})=0
\\
&\forall k\geq 1,\ F(q_{n_k}) =
 p^k_0|\ln |D_{q_{n_0}}||+ p^k_0F(q_{n_0})
+\dots +p^k_{k-1}|\ln |D_{q_{n_{k-1}}}||+ p^k_{k-1}F(q_{n_{k-1}})
\end{split}
\end{equation}

\noindent where the integers $p^k_i$ are given by 
successive euclidean divisions on $dq_{n_k}-\kappa_{max}$:

\begin{equation}\label{divisions-euclidiennes}\begin{split}&dq_{n_k}-\kappa_{max} = p^k_{k-1}dq_{n_{k-1}}+r_{k-1}, \ r_{k-1}<dq_{n_{k-1}},\\
&r_{k-1}= p^k_{k-2}dq_{n_{k-2}}+r_{k-2}, \ r_{k-2}<dq_{n_{k-2}},\\
&\dots, \\
&r_1=p^k_0dq_{n_0}+r^k_0, \ r^k_{0}<dq_{n_{0}}.
\end{split}\end{equation}

\noindent Notice that the $p_i^k$ and $r_0^k$ do depend on $k$.


%


\begin{remk}\label{estim-J_i^k}
Let $i\geq 0$. For all $k\geq i+1$, 

\begin{equation}\label{bornep_i}p_i^k\leq \frac{q_{n_{i+1}}}{q_{n_i}}
\end{equation}

%
%

\noindent Indeed, the integer $p_i^k$ are given by

$$p_{k-1}^k = E\left(\frac{dq_{n_k}-\kappa_{max}}{dq_{n_{k-1}}}\right)$$

\noindent and for $i=0,\dots, k-2$, $p_i^k = E(\frac{r_{i+1}}{dq_{n_i}})$, which satisfies \eqref{bornep_i} by definition of $r_{i+1}$. 
\end{remk}

\begin{lem}\label{fonction-minorante-brjuno}
There exist $C_4>0$ such that the function $F$ satisfies for all $k\geq 0$, 

\begin{equation} 2B(d\alpha)-C_4
\leq \limsup_{k\rightarrow +\infty} \frac{F(q_{n_k})}{q_{n_k}}
\end{equation}

\end{lem}

\dem 
This can be recursively shown. 
Assume that for a fixed $k\geq 1$ and for all $k'\leq k-1$, one has

\begin{equation}\label{minF}\frac{F(q_{n_{k'}})}{q_{n_{k'}}+\frac{\kappa_{max}}{d}}
\geq \sum_{0\leq l\leq k'-1}|\ln |D_{dq_{n_l}}||(\frac{1}{q_{n_l}}-\frac{2}{q_{n_{l+1}}})
\end{equation}

\noindent (which holds for $k=1$). Then

\begin{equation}\begin{split}\frac{F(q_{n_{k}})}{q_{n_{k}}+\frac{\kappa_{max}}{d}}
&\geq 
 \frac{p_{k-1}}{q_{n_{k}}+\frac{\kappa_{max}}{d}}
(F(q_{n_{k-1}})+|\ln |D_{q_{n_{k-1}}}||)\\
&\geq
\frac{p_{k-1}(q_{n_{k-1}}+\frac{\kappa_{max}}{d})}{q_{n_{k}}+\frac{\kappa_{max}}{d}}
\sum_{l\leq k-2}|\ln |D_{dq_{n_l}}||(\frac{1}{q_{n_l}}-\frac{2}{q_{n_{l+1}}})+ \frac{p_{k-1}}{q_{n_{k}}+\frac{\kappa_{max}}{d}}|\ln |D_{q_{n_{k-1}}}||
\\
\end{split}\end{equation}

\noindent Now by definition of $p_{k-1}$,

\begin{equation}\begin{split}
\frac{p_{k-1}(q_{n_{k-1}}+\frac{\kappa_{max}}{d})}{q_{n_{k}}+\frac{\kappa_{max}}{d}}&\geq \frac{(\frac{dq_{n_k}-\kappa_{max}}{dq_{n_{k-1}}}-1)(q_{n_{k-1}}+\frac{\kappa_{max}}{d})}{q_{n_{k}}+\frac{\kappa_{max}}{d}}\\
&=\frac{(q_{n_k}-q_{n_{k-1}}-\frac{\kappa_{max}}{d})(q_{n_{k-1}}+\frac{\kappa_{max}}{d})}{q_{n_{k-1}}(q_{n_k}+\frac{\kappa_{max}}{d})}
\end{split}\end{equation}

\noindent and this quantity 
is greater than 1 since, by assumption on the subsequence $n_k$, 

$$q_{n_k}\geq q_{n_{k-1}}^2 + \zeta q_{n_{k-1}} +\eta$$

\noindent where $\zeta, \eta$ were defined at the beginning of Section \ref{subsequence}. Therefore

\begin{equation}\begin{split}\frac{F(q_{n_{k}})}{q_{n_{k}}+\frac{\kappa_{max}}{d}}&\geq \sum_{l\leq k-2}|\ln |D_{dq_{n_l}}||(\frac{1}{q_{n_l}}-\frac{2}{q_{n_{l+1}}})+|\ln |D_{q_{n_{k-1}}}||\frac{q_{n_k}-q_{n_{k-1}}-\frac{\kappa_{max}}{d}}{q_{n_{k-1}}(q_{n_k}+\frac{\kappa_{max}}{d})}\\
&\geq  \sum_{l\leq k-2}|\ln |D_{dq_{n_l}}||(\frac{1}{q_{n_l}}-\frac{2}{q_{n_{l+1}}})+|\ln |D_{q_{n_{k-1}}}||(\frac{1}{q_{n_{k-1}}}-\frac{2}{q_{n_{k}}})
\end{split}\end{equation}

\noindent (where the last inequality comes from the fact that $q_{n_0}\geq 1+\kappa_{max}(2+\frac{1}{d})$).
Therefore the property \eqref{minF} holds for all $k\geq 0$. 
%
Lemma \ref{grand-div} implies that

$$\frac{F(q_{n_{k}})}{q_{n_{k}}+\frac{\kappa_{max}}{d}}\geq \sum_{0\leq l\leq k}\frac{|\ln |D_{dq_{n_l}}||}{q_{n_l}}-\sum_{0\leq l\leq k}\frac{4\ln (2q_{n_{l+1}}) }{q_{n_{l+1}}}$$

\noindent which implies

$$\frac{F(q_{n_{k}})}{q_{n_{k}}+\frac{\kappa_{max}}{d}}\geq \sum_{0\leq l\leq k}\frac{|\ln |D_{dq_{n_l}}||}{q_{n_l}}-4\sum_{l\geq 0}\frac{\ln q_l }{q_l}-\sum_{l\geq 0} \frac{4\ln 2}{q_l}$$

\noindent Lemma \ref{somme-convergente} then implies 

$$\frac{F(q_{n_{k}})}{q_{n_{k}}+\frac{\kappa_{max}}{d}}\geq \sum_{0\leq l\leq k}\frac{|\ln |D_{dq_{n_l}}||}{q_{n_l}}
-4C_1-4\ln 2 C_0$$


\noindent Finally, by Lemma \ref{somme-extraite},

$$\limsup_{k\rightarrow +\infty}\frac{F(q_{n_{k}})}{q_{n_{k}}}=\limsup_{k\rightarrow +\infty}\frac{F(q_{n_{k}})}{q_{n_{k}}+\frac{\kappa_{max}}{d}}\geq 2B(d\alpha) -2C_2-4C_1-4\ln 2C_0$$

%
%
%
%
%
\noindent where $C_2$ was defined in Lemma \ref{somme-extraite}. Thus one can define $C_4=2C_2
+4C_1+4\ln 2C_0$. $\Box$

\section{An upper bound on the radius of convergence}

In this section, one shall assume the following:

\begin{hyp}\label{argument-coeff}
Assume that there exists $\theta\in\mathbb{R}$ such that for all $k=1,\dots, N$ with $a_k\neq 0$, the complex number $a_k$ has argument $k\theta+\frac{\pi}{2}$. 
\end{hyp}

\noindent
In this case, one can prove a lower bound on the coefficients $\phi_l$ of the linearization in order to bound the radius of convergence from above.

\bigskip
\noindent
One needs the following simple lower bound on all coefficients, including those not corresponding to a small divisor.

\begin{lem}\label{minorgrossire}
For all $r\in\mathbb{N}$, if $|\phi_r|\neq 0$, then $|\phi_r| \geq \min(1,\frac{|a_{\kappa_{min}}|}{4})^{\frac{r}{\kappa_{0}}}$. 

\end{lem}

\dem For all $K$ such that $a_K\neq 0$, it holds that $|\phi_K|\geq \frac{K|a_K|}{4}\geq \frac{|a_{\kappa_{min}}|}{4}$. If $\frac{|a_{\kappa_{min}}|}{4} <1$, then 
$|\phi_K|\geq (\frac{|a_{\kappa_{min}}|}{4})^{\frac{K}{\kappa_0}}$. If $\frac{|a_{\kappa_{min}}|}{4}  \geq 1$ then $|\phi_K|\geq 1$.


\bigskip
\noindent
Let $r\in\mathcal{M},r\geq 2$. Assume that the property holds for all $r'\leq r-1$. Then, either there exists $i\in \{0,\dots,N\}$ such that $r=\kappa_i$, and in this case, 
$|\phi_r| \geq  \frac{a_{\kappa_i}}{4} \geq  \min(1,\frac{|a_{\kappa_{min}}|}{4}) \geq  \min(1,\frac{|a_{\kappa_{min}}|}{4})^{\frac{\kappa_i}{\kappa_0}}$, or there exists $i\in\{0,\dots,N\}$ such that

$$|\phi_r| \geq \frac{|a_{\kappa_i}|}{4}  |\phi_{r-\kappa_i}|\geq \frac{|a_{\kappa_{min}}|}{4}  |\phi_{r-\kappa_i}|
$$

\noindent If 
$\frac{|a_{\kappa_{min}}|}{4}  <1$, then the recurrence assumption implies that
$|\phi_r| \geq (\frac{|a_{\kappa_{min}}|}{4})^{\frac{\kappa_i}{\kappa_0} } \min(1,\frac{|a_{\kappa_{min}}|}{4})^{\frac{r-\kappa_i}{\kappa_0} }\geq  \min(1,\frac{|a_{\kappa_{min}}|}{4})^{\frac{r}{\kappa_{0}}}$. 
If $\frac{|a_{\kappa_{min}}|}{4} \geq 1$ then $|\phi_r| \geq 1$.  $\Box$

\begin{lem} \label{min-prelim}
For all $j> \kappa_0 $ and $p\in\mathbb{N}^*$, there is

$$|\phi_{pj}|\geq \frac{1}{p}(\frac{1}{4})^p|D_j\phi_j|^p.$$

%
%
%

\end{lem}

\dem Notice that

\begin{equation}
|\phi_{pj}|=\frac{1}{|D_{pj}|}|\psi_{pj}|=\frac{1}{|D_{pj}|}\sum_{K\in\{\kappa_0,\dots,\kappa_N\}}|\psi_{K,pj}|
\end{equation}

\noindent (one uses the fact that the $\psi_{K,pj}$ have the same argument for every $K$). Thus, using \eqref{recPsi} with $n=pj,k=(p-1)j$,

\begin{equation}
|\phi_{pj}|\geq \frac{1}{|D_{pj}|}\sum_K \frac{(p-1)j}{pj-K}|\phi_{(p-1)j}||\psi_{K,j}|\geq \frac{(p-1)j}{|D_{pj}| pj}|\phi_{(p-1)j}||\psi_{j}|
\end{equation}

\noindent Iterating this, one obtains

\begin{equation}
|\phi_{pj}|\geq \frac{1}{p|D_{pj}|\dots |D_{2j}|}|\phi_j||\psi_j|^{p-1}= \frac{1}{p|D_{pj}|\dots |D_{2j}|}|\phi_j|^p|D_j|^{p-1}\geq \frac{1}{p}(\frac{1}{4})^p|D_j\phi_j|^p.\ \Box
\end{equation}

%
%
%
%
%
%
%
%
%
%
%
%
%
%

\bigskip
\noindent The following lemma states a better lower bound for the coefficients of the linearization corresponding to a small divisor.

\begin{lem}\label{minor-gros-terme}
For all $k\geq 0$, there is

\begin{equation}\begin{split}|\phi_{dq_{n_{k+1}}}|&\geq  \frac{1}{|D_{dq_{n_{k+1}}}|} \frac{|a_{\kappa_{max}}| }{(k+2)!}
(\frac{1}{4})^{p_0^{k+1}+\dots + p_k^{k+1}}\frac{|D_{dq_{n_0}}|^{p_0^{k+1}}}{p_0^{k+1}}\dots \frac{|D_{dq_{n_{k}}}|^{p_k^{k+1}}}{p_k^{k+1}}\\
&\cdot |\phi_{dq_{n_0}}|^{p_0^{k+1}}\dots |\phi_{dq_{n_k}}|^{p_k^{k+1}} \min(1,\frac{|a_{\kappa_{min}}|}{4})^{\frac{r_0^{k+1}}{\kappa_0}}
\end{split}\end{equation}


\noindent where the integers $p_0^{k+1},\dots, p_k^{k+1},r_0^{k+1}$ were defined in equation \eqref{fonctionF}.
\end{lem}

\dem Taking in the recurrence relation \eqref{rec} the term with $l=q_{n_{k+1}}, K=\kappa_{max},m=1,p_1=k+2,j_1^1 = dp_0^{k+1}q_{n_0} 
,\dots,
j_{k+1}^1 = d
p_k^{k+1}q_{n_k} $ and $j_{k+2}^1=r_0^{k+1}$, in order to have $j_1^1+\dots + j^1_{k+2}=dq_{n_{k+1}}-\kappa_{max}$, one obtains

$$|\phi_{dq_{n_{k+1}}}|\geq \frac{1}{|D_{dq_{n_{k+1}}}|} \frac{\kappa_{max}|a_{\kappa_{max}} |}{(k+2)!} |\phi_{dp_0^{k+1}q_{n_0}}|\dots
|\phi_{dp_k^{k+1}q_{n_k}}||\phi_{r_0^{k+1}}|$$

\noindent
Now apply Lemma \ref{min-prelim} with $j=dq_{n_i}$ and $p=p_i^{k+1}$ to every factor in the right hand side, except the last one, and apply Lemma \ref{minorgrossire} to the last factor. One obtains

\begin{equation}\begin{split}|\phi_{dq_{n_{k+1}}}|&\geq  \frac{1}{|D_{dq_{n_{k+1}}}|} \frac{|a_{\kappa_{max}}| }{(k+2)!}
(\frac{1}{4})^{p_0^{k+1}+\dots + p_k^{k+1}}\frac{|D_{dq_{n_0}}|^{p_0^{k+1}}}{p_0^{k+1}}\dots \frac{|D_{dq_{n_{k}}}|^{p_k^{k+1}}}{p_k^{k+1}}\\
&\cdot |\phi_{dq_{n_0}}|^{p_0^{k+1}}\dots |\phi_{dq_{n_k}}|^{p_k^{k+1}} \min(1,\frac{|a_{\kappa_{min}}|}{4})^{\frac{r_0^{k+1}}{\kappa_0}}.\ \Box
\end{split}\end{equation}

\noindent The following proposition links the radius of convergence of the linearization to the function $F$.

\begin{prop}\label{minoration}
There exists $C\geq 0$ 
such that for all $k$,

$$\frac{1}{dq_{n_{k}}} \ln |D_{dq_{n_k}}\phi_{dq_{n_{k}}} |\geq 
\frac{F(q_{n_k})}{dq_{n_k}}
-C
$$

\noindent Moreover, it is possible to define

$$C=|\ln ( \min(1,\frac{|a_{\kappa_{min}}|}{4}))|+
 \frac{1}{1+N_\mathcal{M}}(2+9  \ln 4)+C_1+[|\ln (\min(1, \frac{ |a_{\kappa_{min}}|}{4}))|+\frac{1}{d}|\ln |a_{\kappa_{max}}||]C_0$$

\end{prop}


\bigskip
\dem
Denote $\tilde{C}=-\ln (\frac{|a_{\kappa_{min}}|}{4})$ if $\frac{|a_{\kappa_{min}}|}{4}<1$ and $\tilde{C}=0$ otherwise. 
Let

\begin{equation}\begin{split}&S_0=0
;\\
&\forall k\geq 1,\  S_k=\sum_{1\leq j\leq k}\frac{(j+1)\ln (j+1)+\ln (q_{n_j})}{q_{n_j}}+\sum_{0\leq j\leq k}\frac{\ln 
4}{q_{n_j}}(1+p_0^j+\dots +p_{j-1}^j)\\&+[\frac{r_0^j}{\kappa_0}|\ln ( \frac{ |a_{\kappa_{min}}|}{4})|+\frac{1}{d}|\ln |a_{\kappa_{max}}||]\sum_{1\leq j\leq k} \frac{1}{q_{n_j}}
\end{split}\end{equation}

%

\noindent Thus, for all $k\geq 0$, $S_k\geq 0$.

\bigskip
\noindent
Let $k\geq 0$. Let us formulate the following recurrence property:

\begin{equation}\label{HR}
\begin{split}
\frac{1}{dq_{n_{k}}} \ln| \phi_{dq_{n_{k}}} |& \geq 
\frac{F(q_{n_k})+|\ln |D_{dq_{n_k}}|}{dq_{n_k}}
- 
S_{k}
-\tilde{C}
\end{split}\end{equation}

\noindent 
First note that this property holds for $k=0$. Indeed, since by Corollary \ref{corsemigp}, $q_{n_0}\geq N_\mathcal{M},$ which implies that $dq_{n_0}\in\mathcal{M}$, then there exists $b_0,b_1,\dots,b_N\geq 0$ such that $dq_{n_0}= \sum_{i=0}^N b_i \kappa_i$. Let $I$ be such that $b_I>0$. Then

$$|\phi_{dq_{n_0}}|\geq \frac{1}{|D_{dq_{n_0}}|} \kappa_I |a_{\kappa_I}| |\phi_{\kappa_I}|^{b_I-1}\prod_{J\neq I} |\phi_{\kappa_J}|^{b_J}$$

\noindent
Now on the other side, for all $0\leq i\leq N$,

$$|\phi_{\kappa_i}|\geq \frac{|a_{\kappa_i}|}{|D_{\kappa_i}|}\geq \frac{|a_{\kappa_{min}}|}{4}
$$

\noindent
hence, if $\frac{|a_{\kappa_{min}}|}{4}<1$,

$$|\phi_{dq_{n_0}}|\geq \frac{1}{|D_{dq_{n_0}}|}  |a_{\kappa_I}| (\frac{|a_{\kappa_{min}}|}{4})^{dq_{n_0}-1}\geq \frac{1}{|D_{dq_{n_0}}|}  (\frac{|a_{\kappa_{min}}|}{4})^{dq_{n_0}}$$

\noindent and if $\frac{|a_{\kappa_{min}}|}{4}\geq 1$,

$$|\phi_{dq_{n_0}}|\geq \frac{1}{|D_{dq_{n_0}}|}  |a_{\kappa_I}|\geq \frac{1}{|D_{dq_{n_0}}|} $$

\noindent
Thus if $\frac{|a_{\kappa_{min}}|}{4}<1$,

$$\frac{1}{dq_{n_0}}\ln |\phi_{dq_{n_0}}|\geq \frac{|\ln |D_{dq_{n_0}}||}{dq_{n_0}}+
\ln (\frac{|a_{\kappa_{min}}|}{4})
\geq \frac{|\ln |D_{dq_{n_0}}||}{dq_{n_0}} - S_0
-\tilde{C}$$

\noindent and if $\frac{|a_{\kappa_{min}}|}{4}\geq 1$,

$$\frac{1}{dq_{n_0}}\ln |\phi_{dq_{n_0}}|\geq \frac{|\ln |D_{dq_{n_0}}||}{dq_{n_0}}
$$

\noindent
therefore the property \eqref{HR} 
holds for $k=0$.

\bigskip
\noindent
Now assume that the recurrence property holds for all $0\leq k'\leq k$, for a fixed $k\geq 0$.
Lemma \ref{minor-gros-terme} implies

\begin{equation}\begin{split}\ln|D_{dq_{n_{k+1}}}\phi_{dq_{n_{k+1}}}| &\geq 
 \ln \left(\frac{|a_{\kappa_{max}}|}{(k+2)!}\right)
-(p_0^{k+1}+\dots +p_k^{k+1})\ln 4-\ln (p_0^{k+1}\dots p_k^{k+1})
\\
&+{p_0^{k+1}}\ln |D_{dq_{n_0}}\phi_{dq_{n_0}}|+\dots+{p_k^{k+1}}\ln |D_{dq_{n_k}}\phi_{dq_{n_k}}|+\frac{r_0^{k+1}}{\kappa_0}\ln \min(1,\frac{|a_{\kappa_{min}}|}{4})\\
\end{split}\end{equation}

\noindent
By recurrence assumption, one infers

\begin{equation}\begin{split}
\frac{1}{dq_{n_{k+1}}} \ln |D_{dq_{n_{k+1}}}\phi_{dq_{n_{k+1}}}| & \geq 
-\frac{(k+2)\ln (k+2)}{dq_{n_{k+1}}}+\frac{\ln |a_{\kappa_{max}}|}{dq_{n_{k+1}}}-\frac{p_0^{k+1}+\dots + p_k^{k+1}}{dq_{n_{k+1}}}\ln 4
\\
 &+\frac{p_0^{k+1} dq_{n_0}}{dq_{n_{k+1}}}(\frac{F(q_{n_0})+|\ln |D_{q_{n_0}}||}{dq_{n_0}}-S_0-\tilde{C})+\dots \\
 &+ \frac{p_k^{k+1}dq_{n_k} }{dq_{n_{k+1}}}(\frac{F(q_{n_k})+|\ln |D_{q_{n_k}}||}{dq_{n_k}}-S_{k}-\tilde{C})\\
 &+\frac{r_0^{k+1}}{d\kappa_0q_{n_{k+1}}}\ln \min(1,\frac{|a_{\kappa_{min}}|}{4})-\frac{\ln (p_0^{k+1}\dots p_k^{k+1})}{dq_{n_{k+1}}}
\end{split}\end{equation}

\noindent thus by definition of $F$,

\begin{equation}\begin{split}
\frac{1}{dq_{n_{k+1}}} \ln |D_{dq_{n_{k+1}}}\phi_{dq_{n_{k+1}}}| 
&\geq \frac{F(q_{n_{k+1}})}{dq_{n_{k+1}}} -\frac{(k+2)\ln (k+2)}{dq_{n_{k+1}}}+\frac{\ln |a_{\kappa_{max}}|}{dq_{n_{k+1}}}-\frac{p_0^{k+1}+\dots + p_k^{k+1}}{dq_{n_{k+1}}}\ln 4
\\
& +  \frac{r_0^{k+1}}{d\kappa_0q_{n_{k+1}}}\ln \min(1,\frac{|a_{\kappa_{min}}|}{4})-\tilde{C} - \frac{p_0^{k+1}dq_{n_0}}{dq_{n_{k+1}}}S_0 -\dots - \frac{p_k^{k+1}dq_{n_k}}{dq_{n_{k+1}}}S_k\\
&-\frac{\ln (p_0^{k+1}\dots p_k^{k+1})}{dq_{n_{k+1}}}\\
\end{split}\end{equation}

\noindent Now for all $i=0,\dots, k$, $p_i\leq \frac{q_{n_{i+1}}}{q_{n_i}}$ therefore

$$\frac{\ln (p_0^{k+1}\dots p_k^{k+1})}{dq_{n_{k+1}}}\leq \frac{\ln (\frac{q_{n_{k+1}}}{q_{n_0}})}{dq_{n_{k+1}}}$$

\noindent and

$$\frac{1}{dq_{n_{k+1}}} \ln |D_{dq_{n_{k+1}}}\phi_{dq_{n_{k+1}}}| \geq  \frac{F(q_{n_{k+1}})}{dq_{n_{k+1}}}-\tilde{C}-S_{k+1}$$

\noindent 
Therefore, the property \eqref{HR}
holds for all $k\geq 0$. 

\bigskip
\noindent
Finally note that the partial sums
$S_k$ converge, since from one side, by choice of the subsequence $q_{n_j}$ and corollary \ref{somme-suite-extraite},

$$\sum_{1\leq j\leq k}  \frac{(j+1) \ln (j+1)+\ln (q_{n_j})}{q_{n_{j}}}\leq \sum_{1\leq j\leq k}  \frac{1}{\sqrt{q_{n_{j}}}} +C_1\leq \sum_{1\leq j\leq k}\frac{1}{q_{n_{j-1}}}+C_1
\leq \frac{1}{1+N_\mathcal{M}}+C_1
$$

\noindent and from the other side, by Remark \ref{estim-J_i^k},

\begin{equation}\begin{split}
&
\sum_{j=1}^k \frac{\ln 4}{q_{n_j}}(1+p_0^j+\dots + p_{j-1}^j)\\
&\leq 
\sum_{j=1}^k \frac{\ln 4}{q_{n_j}}(1+\frac{q_{n_1}}{q_{n_0}}+\dots + \frac{q_{n_{j-1}}}{q_{n_{j-2}}}+\frac{q_{n_j}}{q_{n_{j-1}}})\\
&\leq \ln 4(
 \sum_{j=1}^k \frac{(j+1)q_{n_{j-1}}}{q_{n_j}}+\frac{1}{q_{n_{j-1}}})\\
&\leq \ln 4(
\sum_{j=1}^k \frac{j+2}{q_{n_{j-1}}})\\
&\leq \ln 4(\frac{8}{q_{n_0}} + \sum_{j=2}^{k-1}\frac{j+3}{q_{n_j}})\leq  \ln 4(\frac{8}{q_{n_0}} 
+ \sum_{j=2}^{k-1}\frac{1}{q_{n_{j-1}}})\\
&\leq \frac{9\ln 4}{1+N_\mathcal{M}}
\end{split}\end{equation}

\noindent
(this last inequality used Corollary \ref{somme-suite-extraite}). Also, 

$$\frac{r_0^j}{\kappa_0}|\ln ( \frac{ |a_{\kappa_{min}}|}{4})|\sum_{1\leq j\leq k} \frac{1}{q_{n_j}}\leq |\ln ( \frac{ |a_{\kappa_{min}}|}{4})|\sum_{1\leq j\leq k} \frac{q_{n_0}}{q_{n_j}}\leq |\ln ( \frac{ |a_{\kappa_{min}}|}{4})|\sum_{j\geq 0}\frac{1}{q_{n_j}}\leq C_0|\ln ( \frac{ |a_{\kappa_{min}}|}{4})|$$

\noindent 
Thus, let 

$$C=\tilde{C}+ \frac{1}{1+N_\mathcal{M}}(2+9  \ln 4)+C_1+[|\ln ( \frac{ |a_{\kappa_{min}}|}{4})|+\frac{1}{d}|\ln |a_{\kappa_{max}}||]C_0$$

\noindent then $C$ satisfies the statement of this proposition.
$\Box$

\begin{thm} The radius of convergence of $\Phi$ is bounded from above by $\exp(-\limsup_{k\rightarrow +\infty}\frac{F(q_{n_k})}{dq_{n_k}}+C)$ where $C\geq 0$ was defined in Proposition \ref{minoration}.
It is also bounded from above by $\exp(-\frac{2}{d}B(d\alpha) + C_\mathcal{M})$ where $C_\mathcal{M} = \frac{C_4}{d}+C$, with $C_4$ defined in Lemma \ref{fonction-minorante-brjuno}.

%
%
%

\end{thm}

\dem Let $\rho$ be the radius of convergence of $\Phi$, then $\rho^{-1}= \limsup_{j\rightarrow +\infty} |\phi_j|^{\frac{1}{j}}$ therefore

$$-\ln \rho \geq  \limsup_{k\rightarrow +\infty} \frac{1}{dq_{k}}\ln |\phi_{dq_{k}}|\geq  \limsup_{k\rightarrow +\infty} \frac{1}{dq_{n_k}}\ln |\phi_{dq_{n_k}}|
\geq \limsup_{k\rightarrow +\infty} \frac{|\ln |D_{q_{n_k}}||}{dq_{n_k}}+\frac{F(q_{n_k})}{dq_{n_k}}-C$$ 

\noindent where $C$ was defined in Proposition \ref{minoration}. 
By Lemma \ref{fonction-minorante-brjuno}, 

$$-\ln\rho \geq \frac{2}{d}B(d\alpha)
-\frac{C_4}{d}-C.$$

One can then define $C_\mathcal{M}= \frac{C_4}{d}+C
$.
 $\Box$

\section{A lower bound on the radius of convergence}

In this section, the second part of the main result is proved.
The assumption \ref{argument-coeff} on the coefficients $a_K$ is relaxed.

\begin{thm}\label{thm-minor-rayon}
The radius of convergence is at least $\exp(-C'-
\sum_{l\geq 0}  \frac{2\ln q_{{l+1}}}{dq_{l}})$, where $C'\geq 0$ is defined by

$$C'= 
\ln (|a_{\kappa_{max}}|) + r+C_0$$

\noindent if $|a_{\kappa_{max}}|>1$ and

$$C'= 
r+C_0$$

 \noindent otherwise, with $r>0$ only depending on the Fourier modes of the trigonometric polynomial $A$ and $C_0$ defined in Lemma \ref{somme-convergente}.
\end{thm}

\dem
Let $w$ be an analytic solution of the functional equation

$$w(z) = \sum_{k\in\{\kappa_0,\dots,\kappa_N\}}(ze^{w(z)})^k$$

\noindent and let $R$ its radius of convergence.
Expanding $w$ in its Taylor series, $w(z) = \sum_{n\geq 0} \sigma_n z^n$, one obtains the following relation between the coefficients $\sigma_n$ (it is the same relation as between the coefficients $\phi_j$, only replacing $a_k$ by $i^{k+1}$ and without small divisors):

\begin{equation}
\sigma_l = \sum_{\substack{K\in\\
\kappa_0,\dots, \\
\kappa_N\},\\
K\leq l}}  [\delta_{K,l}+\sum_{\substack{m=1,\dots, K\\
m\leq l-K}} C^m_K \sum_{p_1=1}^{l-K} \dots \sum_{p_m= 1}^{l-K}\frac{1}{p_1!\dots p_m!}\sum_{
\substack{j_1^1,\dots , j^1_{p_1},\\
\vdots\\
j^m_1,\dots, j^m_{p_m},\\
j_1^1+\dots + j^1_{p_1} +\dots\\
 + j^m_1 +\dots 
+ j^m_{p_m} = l-K}} \sigma_{j^1_1}\dots \sigma_{j^1_{p_1}}\dots \sigma_{j^m_1}\dots \sigma_{j^m_{p_m}}
\big]
\end{equation} 

\noindent
One can recursively show that the $\sigma_n$ are non negative real numbers. Moreover the function $w$ is analytic, therefore $\limsup_{n\rightarrow +\infty} \frac{-1}{n}\ln \sigma_n$ is equal to $-\ln R$.

%
%
%
%
%
%
%
%
%
%
%

\noindent 
The function $g$ giving the upper bound is defined as follows: for all $k\geq 0$, let $g_k$ be the function defined by Davie's lemma \ref{Davie2}. Then,
for all integer $\kappa_0\leq j <q_0$, let 

$$g(j)=
 j\ln ( |a_{\kappa_{max}}|) +36j
$$

\noindent  if $|a_{\kappa_{max}}|>1$, and 

$$g(j)=
36j
$$

\noindent otherwise.
Now let $j\geq q_0$ and 
assume that $k$ is the greatest index such that $q_{k}\leq j$; let

$$g(j)=
\sum_{l=0}^k  2g_l(j)\ln q_{l+1}+ j\ln ( |a_{\kappa_{max}}|) +36j
$$

\noindent  if $|a_{\kappa_{max}}|>1$, and 

$$g(j)=
\sum_{l=0}^k  2g_l(j)\ln q_{l+1} +36j
$$

\noindent otherwise.
%
%
%
%
%
%
The function $g$ is increasing.
Moreover for all $j_1,j_2\geq \kappa_0$, as a consequence of Davie's lemma,

\begin{equation}\label{proprit-g}
g(j_1)+ g(j_2)\leq 
g(j_1+ j_2)
\end{equation}



%
%
%
%
%
%

\bigskip
\noindent
Now let us prove that $|\phi_j|\leq \sigma_j e^{g(j)}$ for all $j\geq \kappa_0$. 
First, 
if $|D_{\kappa_0}|^{-1}|a_{\kappa_{max}}|>1$, then

$$|\phi_{\kappa_0}|=|D_{\kappa_0}|^{-1}|a_{\kappa_0}|
\leq |D_{\kappa_0}|^{-1}|a_{\kappa_{max}}|
\leq e^{g(\kappa_0)}$$

\noindent and $|\phi_{\kappa_0}|\leq 1=e^{g(\kappa_0)}$ otherwise.
Assume that this holds for all $j'\leq j-1$ and consider $|\phi_{j}|$. 
The relation \eqref{rec} implies

\begin{equation}\label{crochet}
|\phi_{j}| \leq |D_j^{-1} |\sum_{\substack{K\in\\
\kappa_0,\dots, \\
\kappa_N\},\\
K\leq j}}
| a_K|[\delta_{K,j}+
\sum_{\substack{m=1,\dots, K\\
m\leq j-K}} C^m_K 
\sum_{p_1=1}^{j-K} \dots \sum_{p_m= 1}^{j-K}\frac{1}{p_1!\dots p_m!}
\sum_{
\substack{j_1^1,\dots , j^1_{p_1},\\
\vdots\\
j^m_1,\dots, j^m_{p_m},\\
j_1^1+\dots + j^1_{p_1} +\dots\\
 + j^m_1 +\dots 
+ j^m_{p_m} = j-K}} |\phi_{j^1_1}|\dots |\phi_{j^1_{p_1}}|\dots |\phi_{j^m_1}|\dots |\phi_{j^m_{p_m}}|
]
\end{equation} 

\noindent
hence, by recurrence assumption,

\begin{equation}
|\phi_{j}| \leq |D_j^{-1}|
 \sum_{\substack{K\in\\
\kappa_0,\dots, \\
\kappa_N\},\\
K\leq j}}
| a_K|[\delta_{K,j}+
\sum_{\substack{m=1,\dots, K\\
m\leq j-K}} C^m_K 
\sum_{p_1=1}^{j-K} \dots \sum_{p_m= 1}^{j-K}\frac{1}{p_1!\dots p_m!}
\sum_{
\substack{j_1^1,\dots , j^1_{p_1},\\
\vdots\\
j^m_1,\dots, j^m_{p_m},\\
j_1^1+\dots + j^1_{p_1} +\dots\\
 + j^m_1 +\dots
+ j^m_{p_m} = j-K}} \sigma_{j^1_1}\dots 
\sigma_{j^m_{p_m}}e^{g(j_1^1)+\dots g(j_{p_m}^m)}
]
\end{equation} 

\noindent 
therefore

$$|\phi_{j}| \leq |D_j^{-1}| |a_{\kappa_{max}}| 
e^{g(j-\kappa_0)}
\sigma_j$$

\noindent
We shall distinguish two cases:

\begin{itemize}


\item if $j\geq \kappa_0$ and $|D_j|^{\frac{1}{2}}\geq \frac{1}{6}$, then

$$ |D_j^{-1}| |a_{\kappa_{max}}| 
e^{g(j-\kappa_0)}
\sigma_j\leq 36|a_{\kappa_{max}}| 
e^{g(j-\kappa_0)}
\sigma_j
\leq e^{g(j)}
\sigma_j
$$

\item otherwise there exists $k\geq 0$ such that $\frac{1}{6q_{{k+1}}}\leq |D_j| ^{\frac{1}{2}}<\frac{1}{6q_{k}}$, and in this case,
$|D_j^{-1}|\leq 36q_{{k+1}}^2$. Moreover, $j\in\mathcal{M} $ implies that $d$ divides $j$; 
then $j\in A_k(j)$ therefore, by construction of $g_k$,

$$g_k(j)=g_k(j-\kappa_0)+1$$

\noindent therefore 

$$\ln |D_j^{-1}|+\ln |a_{\kappa_{max}}|+g(j-\kappa_0)\leq 2\ln q_{k+1}+\ln 36 +\ln |a_{\kappa_{max}}|+g(j-\kappa_0)\leq g(j)$$


\noindent thus

$$|\phi_j|\leq  |D_j^{-1}| |a_{\kappa_{max}}| 
e^{g(j-\kappa_0)}\sigma_j\leq e^{g(j)}\sigma_j $$

%
%
%
%
%
%
%
%
%
%
%
%
%
%

\end{itemize}

\noindent The recurrence is finished. 
Thus, for all $j\geq \kappa_0$,

\begin{equation}\begin{split}\frac{1}{j}\ln |\phi_j| &\leq \frac{1}{j}g(j)+\frac{1}{j}\ln \sigma_j 
\end{split}\end{equation}

\noindent Now 

$$\frac{1}{j}g(j)\leq 
\sum_{0\leq l\leq k}  \frac{2}{j}g_l(j)\ln q_{l+1}
+\delta \ln ( |a_{\kappa_{max}}|)+36$$

\noindent where $\delta=1$ if $|a_{\kappa_{max}}|>1$ and $0$ otherwise. Therefore by Lemma \ref{Davie2},

$$\frac{1}{j}g(j)
\leq \sum_{0\leq l\leq k}  (\frac{2}{dq_l} + \frac{16}{dq_{l+1}})\ln q_{l+1}
+\delta \ln ( |a_{\kappa_{max}}|)+36
$$

\noindent 
Thus,

\begin{equation}\begin{split}
\limsup_{j\rightarrow +\infty} \frac{1}{j}\ln |\phi_j| &\leq
 \limsup_{j\rightarrow +\infty} \frac{1}{j}g(j)+\limsup_{j\rightarrow +\infty}\frac{1}{j}\ln \sigma_j \\
& \leq 
\sum_{l\geq 0}  \frac{2\ln (q_{{l+1}})}{dq_{l}}+\sum_{l\geq 0}  \frac{16\ln (q_{{l+1}})}{dq_{l+1}}+\delta \ln ( |a_{\kappa_{max}}|)+36
+\limsup_{j\rightarrow +\infty}\frac{1}{j}\ln \sigma_j \\
 \end{split}\end{equation}

%
%
%

\noindent
Thus,

$$\limsup_{j\rightarrow +\infty} \frac{1}{j}\ln |\phi_j| \leq 
2\sum_{l\geq 0}  \frac{\ln q_{{l+1}}}{dq_{l}}
+ C'$$

\noindent where

$$C'= 
 \ln(|a_{\kappa_{max}}|) + \ln R+\frac{16 C_1}{d}+36$$

\noindent if $|a_{\kappa_{max}}|>1$ and $R>1$,

$$C'= 
 \ln R+\frac{16 C_1}{d}+36$$

\noindent if $R>1 $ and $|a_{\kappa_{max}}|\leq 1$, and

$$C'=\frac{16 C_1}{d}+36$$

\noindent otherwise. It only remains to use the characterization of the radius as a function of $\limsup_{j\rightarrow +\infty} \frac{1}{j}\ln |\phi_j|$.
 $\Box$

%
%
%
%
%
%
%
%
%
%
%

\section{Appendix}


\noindent
Here we give a proof of equation \eqref{equation-Phi} in two lemmas.

\begin{lem}\label{2ecomposante}
Let $H(z,\lambda)=(h(z,\lambda),h_2(z,\lambda))$ be the linearization. Then $h_2(z,\lambda)=\frac{h(z,\lambda)}{h(\lambda^{-1}z,\lambda)}$.
\end{lem}

\dem Let $F(z,\lambda)=(\lambda z \prod_{k=1}^N e^{ia_k z^k}, \lambda \prod_{k=1}^N e^{ i a_k z^k})$ be the function generating the system, and $R(z,\lambda)=(\lambda z,\lambda)$. Note that if $f(z,\lambda)=\prod_{k=1}^N e^{ i a_k z^k}$ then $F(z,\lambda) = (z f(z,\lambda) , f(z,\lambda))$. We shall expand the identity $H = F\circ H\circ R^{-1}$ to infer the desired identity.
On one side,

\begin{equation}\begin{split}
F\circ H\circ R^{-1} (z,\lambda)& = F\circ H(\lambda^{-1} z,\lambda) = F(h(\lambda^{-1} z,\lambda), h_2(\lambda^{-1} z,\lambda)) \\
&= ( h(\lambda^{-1} z,\lambda) f(h(\lambda^{-1}z,\lambda),h_2(\lambda^{-1} z,\lambda)), 
 f(h(\lambda^{-1}z,\lambda),h_2(\lambda^{-1} z,\lambda)))
 \end{split}\end{equation}
 
 \noindent
 By matching the components of $F\circ H\circ R^{-1} $ with those of $H$, one infers that
 
 $$h(z,\lambda) = h(\lambda^{-1} z,\lambda) f(h(\lambda^{-1}z,\lambda),h_2(\lambda^{-1} z,\lambda))$$
 
 \noindent and
 
 $$h_2(z,\lambda) = f(h(\lambda^{-1}z,\lambda),h_2(\lambda^{-1} z,\lambda)))$$
 
 \noindent Those two identities imply that for all $(z,\lambda)$,
 
 $$h_2(z,\lambda) = \frac{h(z,\lambda)}{h(\lambda^{-1} z,\lambda) }. \Box$$

\begin{lem}
Let $h(z,\lambda)=ize^{\Phi_\lambda(z)}$, then equation \eqref{equation-Phi} holds.
\end{lem}

\dem Let $F(z,\lambda)=(\lambda z \prod_{k=1}^N e^{ia_k z^k}, \lambda \prod_{k=1}^N e^{ i a_k z^k})$ be the function generating the system, and let $R(z,\lambda)=(\lambda z,\lambda)$. 
By definition, $F\circ H=H\circ R$. Now using Lemma \ref{2ecomposante}, 

$$F\circ H(z,\lambda) = (\frac{h^2(z,\lambda)}{h(\lambda^{-1} z,\lambda)}\prod_{k=1}^N e^{ia_k h(z,\lambda)^k},
\frac{h(z,\lambda)}{h(\lambda^{-1} z,\lambda)}\prod_{k=1}^N e^{ia_k h(z,\lambda)^k})$$

\noindent and $H\circ R(z,\lambda) = (h(\lambda z,\lambda), \frac{h(\lambda z,\lambda)}{h(z,\lambda)})$. By matching the components and taking the logarithm modulo $2 i \pi$, one obtains equation \eqref{equation-Phi}. $\Box$

%
%
%
%
%

\end{document}